\documentclass[reqno]{amsart}

\usepackage{enumerate}
\usepackage{amssymb}

\usepackage[colorlinks=true]{hyperref}
\hypersetup{urlcolor=blue, citecolor=cyan}
\hypersetup{citecolor=blue}

\numberwithin{equation}{section}

\newtheorem{thm}{Theorem}[section]
\newtheorem{lem}[thm]{Lemma}
\newtheorem{prop}[thm]{Proposition}

\theoremstyle{definition}
\newtheorem{rem}[thm]{Remark}

\newcommand\R{{\mathbb R}}
\newcommand\Rn{{\mathbb R}^N}
\newcommand\C{{\mathbb C}}
\newcommand\N{{\mathbb N}}

\newcommand\MScN[1]{\href{http://www.ams.org/mathscinet-getitem?mr=#1}{\nolinkurl{(#1)}}}
\newcommand\DOI[1]{\href{http://dx.doi.org/#1}{(doi: \nolinkurl{#1})}}
\newcommand\LINK[1]{\href{#1}{(link: \nolinkurl{#1})}}

\title{ Existence of standing waves for the complex Ginzburg-Landau equation}

\author[Rolci Cipolatti, Fl\'avio Dickstein and Jean-Pierre Puel]{}



\thanks{Fl\'avio Dickstein  was partially supported by CNPq (Brasil).}
\thanks{This work has been done while Jean-Pierre Puel was visiting Universidade Federal do Rio de Janeiro as a "Professor Visitante Especial" of the "Programa Ci\^encia sem Fronteiras" of Capes/CNPq (Brasil).}

\begin{document}
\maketitle

\centerline{\scshape Rolci Cipolatti}
\medskip
{\footnotesize
 \centerline{Instituto de Matem\'atica}
 \centerline{Universidade Federal do Rio de Janeiro}
   \centerline{ Caixa Postal 68530}
   \centerline{21944--970 Rio de Janeiro, R.J., Brazil}
   \centerline{email address: {\href{mailto:cipolatti@ufrj.br}{\tt cipolatti@ufrj.br}}}
}

\medskip

\centerline{\scshape Fl\'avio Dickstein}
\medskip
{\footnotesize
 \centerline{Instituto de Matem\'atica}
 \centerline{Universidade Federal do Rio de Janeiro}
   \centerline{ Caixa Postal 68530}
   \centerline{21944--970 Rio de Janeiro, R.J., Brazil}
   \centerline{email address: {\href{mailto:flavio@labma.ufrj.br}{\tt flavio@labma.ufrj.br}}}
  %  \centerline{URL:  {\href{http://dma.im.ufrj.br/docentes/flavio/}{\tt http://dma.im.ufrj.br/docentes/flavio/}} }
}

\medskip

\centerline{\scshape Jean-Pierre~Puel}
\medskip
{\footnotesize
 \centerline{Universit\'e de Versailles Saint-Quentin}
 \centerline{LMV, CNRS UMR 8100}
   \centerline{45 avenue des Etats Unis}
   \centerline{78035 Versailles, France}
   \centerline{email address: {\href{mailto:jppuel@math.uvsq.fr}{\tt jppuel@math.uvsq.fr}}}
 %   \centerline{URL:  {\href{http://www.math.univ-paris13.fr/ann/indiv/index.php?clef=WFoySGl8Tr}{\tt http://www.math.univ-paris13.fr/ann/}} }
}

\begin{abstract}
We prove the existence of non-trivial standing wave solutions of the complex Ginzburg-Landau equation $\phi _t - e^{ i\theta }(\rho I- \Delta) \varphi - e^{ i\gamma  }  |\phi |^\alpha \varphi =0 $ in $\Rn$, where $(N-2)\alpha <4$, $\theta ,\gamma \in (-\pi /2,\pi /2)$ and $\rho >0$.  Analogous result is obtained in a ball $\Omega \in\Rn$ for $\rho >-\lambda _1$, where $\lambda _1$ is the first eigenvalue of the Laplace operator with Dirichlet boundary conditions. 
\end{abstract}%\subjclass[2010] {35Q56, 35C08}

% \keywords{Standing waves, complex Ginzburg-Landau equation}

\section{Introduction}
The complex Ginzburg-Landau equation 
\begin{equation} \label{GCGL} 
\psi _t =z_1 \Delta \psi + z_2  |\psi |^\alpha \psi + z_3\psi ,
\end{equation} 
for $\alpha =2$, $z_1, z_2,  z_3\in \C$, with $\Re z_1\ge 0$ was proposed independently by DiPrima, Eckhaus, Segel \cite{des} and  Stewartson, Stuart \cite{StewartsonS} to model the  interaction of plane waves in fluid flows and plays a central role in the study of the development of nonlinear instabilities in fluid dynamics. See~\cite{DoeringGHN, CrossH, vanSaarloosH} and the references cited therein for a discussion of various problems where the complex Ginzburg-Landau equation applies. Local (global for $\Re z_2<0$) well-posedness of~\eqref{GCGL} (for $\alpha >0)$ were derived in both $\R^N $ and a domain $\Omega \subset \R^N $, under various boundary conditions and assumptions on the parameters, in~\cite{DoeringGL, GinibreVu, GinibreVd, LevermoreO, LevermoreOd, MischaikowM, OkazawaYu, OkazawaYd, OkazawaYt, OkazawaYq}. 

The existence of special solutions  of \eqref{GCGL}  (holes, fronts, pulses, sources, sinks, etc) is discussed in numerous works, see e.g.~\cite{ChungC, Cruz-PachecoLL, DescalziAT, Doelman, LanGC, LegaF, MancasC, MohamadouNK, PoppSAK, PoppSK, vanSaarloosH}. We look for standing wave solutions.  Replacing $\varphi $ by $e^{i\eta t} \varphi $ for some $\eta\in \R$ and rescaling the equation,  we rewrite \eqref{GCGL}  as  
\begin{equation} \label{CGL} 
\partial _t \varphi +e^{i\theta }(\rho \varphi -\Delta \varphi)=e^{i\gamma }|\varphi |^\alpha\varphi ,
\end{equation} 
where $\rho \in \R$. 
Given $ \omega \in \R$, a standing wave of the form $\varphi = e^{i\omega t} u (x)$ is a solution of  \eqref{CGL} if and only if $u$ satisfies 
\begin{equation} \label{sw} 
i\omega  u +e^{i\theta }(\rho u -\Delta u)=e^{i\gamma }|u |^\alpha u .
\end{equation}
Plane waves $\varphi =e^{i(kx-\omega t)}$, where $k,\omega \in \R$ are particular standing waves. It is easy to see that \eqref{CGL} admits plane wave solutions in $\Rn$ for all values of $\rho $, $\theta $, $\gamma $ and $\alpha $. Stationary solutions are also standing waves of special kind. If $\omega =0$ and $u\ne 0$ then necessarily $\sin \gamma =\sin\theta $, so that  equation~\eqref{sw} reduces to the nonlinear elliptic equation $\rho u-\Delta u=\pm |u|^\alpha u$. The case of the nonlinear Schr\"o\-din\-ger equation $\theta= \pm \gamma  =\pm \frac {\pi } {2}$ leads to the equation $\Delta u +  |u|^\alpha u -(\rho \mp \omega )u=0$.   

We will obtain solutions that are different from these particular ones. In fact, using well known results of the theory of nonlinear elliptic equations for the case $\omega =0$ and $\theta =\gamma $, we show the existence of nontrivial standing wave solutions for $\theta \ne \gamma $ by a perturbation argument, as we describe below.     
  
Equation \eqref{sw} will be considered both in the whole space $\Omega =\Rn$ or in a ball $\Omega \in \Rn$ with Dirichlet boundary condition, for $N\ge 1$. We suppose $\theta ,\gamma \in (-\pi /2,\pi /2)$ and $\alpha $ subcritical, i.e.  
\begin{equation} \label{sub_crit} 
   (N-2)\alpha <4,
\end{equation} 
which includes the relevant case $\alpha =2$, for $N\le 3$. For $\theta =\gamma $ and $\omega =0$, \eqref{sw} reduces to 
 \begin{equation} \label{ellip} 
\rho u-\Delta u-  |u |^\alpha u =0 .
\end{equation} 

Consider first $\Omega =\Rn$, in which case we assume that $\rho >0$. It is then known that  \eqref{ellip} has a unique positive radially symmetric solution $U\in C^2 (\Rn )\cap C_0(\Rn)$. In fact, $U\in H_{\rm rad}^2 (\Rn )$, the subspace of  radial functions of $H^2 (\Rn )$.  Note that \eqref{ellip} is phase invariant, i.e., $Ue^{i\beta }$ is also a solution for all $\beta \in\R$.  We prove the following result, in which the Hilbert spaces are real, but composed of complex-valued functions.
\begin{thm} \label{thm1} 
Assume \eqref{sub_crit} holds and suppose $\rho >0$. Let $U\in H_{\rm rad}^2 (\Rn )$ be the unique positive radial solution of \eqref{ellip}. Given 
$\theta \in (-\pi /2,\pi /2)$ and $\beta\in\R$ there exists $0<\varepsilon<\min\{\pi/2-\theta,\pi /2+\theta \}$ and a $C^1$ mapping $g:(\theta -\varepsilon,\theta +\varepsilon)\to \R\times H_{\rm rad}^2(\Rn )$, $g(\gamma )=(\omega_\gamma  ,u_\gamma )$, satisfying $\omega_\theta =0$, $u_\theta =Ue^{i\beta }$ and such that $\varphi _\gamma =e^{i\omega _\gamma t}u_\gamma $ is a solution of \eqref{CGL}. 
\end{thm} 
In the bounded domain case of the unitary ball $\Omega $ of $\Rn$, we suppose that 
\begin{equation} \label{rho2} 
\rho>-\lambda _1,
\end{equation} 
where $\lambda _1$ is the first eigenvalue associate to the Laplace-Dirichlet operator in $\Omega $. As in the case of the whole space, \eqref{ellip} admits a unique  positive solution $U\in H^2(\Omega )\cap H^1_0(\Omega )$, which is radial and radially decreasing.  The following result is analogous to Theorem~\ref{thm1}.  
\begin{thm} \label{thm2} 
Assume \eqref{sub_crit}, \eqref{rho2} hold and let $U\in H^2(\Omega )\cap H^1_0(\Omega ))$ be the positive solution of \eqref{ellip}. Given 
$\theta \in (-\pi /2,\pi /2)$ and $\beta\in\R$ there exists $0<\varepsilon<\min\{\pi/2-\theta,\pi /2+\theta \}$ and a $C^1$ mapping $g:(\theta -\varepsilon,\theta +\varepsilon)\to \R\times (H^2(\Omega )\cap H^1_0(\Omega ))$, $g(\gamma )=(\omega_\gamma  ,u_\gamma )$, satisfying $\omega_\theta =0$, $u_\theta =Ue^{i\beta }$ and such that $\varphi _\gamma =e^{i\omega _\gamma t}u_\gamma $ is a solution of \eqref{CGL}. 
\end{thm} 
In the proofs of Theorem~\ref{thm1} and Theorem~\ref{thm2} we apply the Implicit Function Theorem to $F(w,u,\gamma )= i\omega  u +e^{i\theta }(\rho u -\Delta u)-e^{i\gamma }|u |^\alpha u$ in a neighborhood of $w=0$, $u=Ue^{i\beta }$ and $\gamma =\theta $. Analogous approach was considered in \cite{cdw} to obtain standing wave solutions to \eqref{CGL} in a bounded domain for $\alpha $ small,  where an eigenvector of the Laplace-Dirichlet operator is used as a starting point.  Our point of view allows us to obtain solutions for $\alpha $ satisfying \eqref{sub_crit}  and for the case of the whole space. We are lead to study the linearized operator $L_\beta =\partial _u F(0,Ue^{i\beta },\theta )$ in an appropriate setting. In fact,  it will be sufficent to consider $L =\partial _u F(0,U,\theta )$, see Section~\ref{proofs}. 

We address some comments about the hypothesis in Theorem~\ref{thm1} and Theorem~\ref{thm2}. 
\begin{rem} 

\begin{enumerate}
\item The assumption $\theta \in (-\pi /2,\pi /2)$ yields an accretive linear operator associated to the problem and corresponds to $\Re z_1 >0$ in \eqref{GCGL}. We also obtain $\gamma \in (-\pi /2,\pi /2)$, i.e., standing waves appear in the focusing case. In the defocusing case $\gamma \in (\pi /2,3\pi /2)$, multiplying the equation by $ \overline{\varphi }$ and integrating, we see that $\|\varphi (t)\|_{L^2 (\Rn ) }$ decreases in time.  Thus there cannot be any non-trivial standing wave in that case.
\item The restriction to radial solutions in Theorem~\ref{thm1} seems to be necessary in our proof. It ensures the  compactness of the linear operator $K$ introduced in the proof. It also ensures that $\ker L$ is one-dimensional, which allows for the application of the Implicit Function Theorem. As discussed in Section~\ref{Rn},  $\ker L$ is $(N+1)$-dimensional in $L^2 (\Rn )$. 
\item The assumption that $\Omega $ is a ball in Theorem~\ref{thm2} ensures that $\ker L$ is one-dimensional. We don't know if this is true in general. The standing waves in Theorem~\ref{thm2} can be constructed such that they are radially symmetric, see Remark~\ref{rem}. 
\end{enumerate}       
\end{rem} 

This paper is organized as follows. In Section~\ref{start} we recall some well stablished properties of the positive solution $U$, both in the bounded and in the unbounded domain cases. A spectral analysis of the operator $L $  is developed in Section~\ref{ball} for the case where $\Omega $ is a ball, and in Section~\ref{Rn} when $\Omega $ is  the whole space. Finally, in Section~\ref{proofs} we prove Theorem~\ref{thm1} and Theorem~\ref{thm2}.   

\section{The starting point $\theta =\gamma $}\label{start} 
In this section we recall some well known properties of solutions $u\in H^1_0 (\Omega ) $ of \eqref{ellip} which will be useful later. 

We consider first the case where $\Omega $ is a ball and we assume \eqref{rho2}.  Then  \eqref{ellip} admits infinitely many real solutions and, in particular, a positive radially symmetric solution $U$ \cite{ar}. Equation \eqref{ellip} is phase invariant: if $u$ solves \eqref{ellip} so does $e^{i\beta }u$ for all $\beta \in \R$.  Non-radial complex solutions in $\Rn$ were obtained in \cite{lions}. 

The positive solution $U$, which was shown to be unique in  \cite{kwong}, can be obtained by ode's methods \cite{cazenave}. It can also be derived by solving the minimization problem
\begin{equation} \label{min} 
\min_{u\in S} \int \rho |u|^{2}+ |\nabla u|^{2},
\end{equation}  
where 
\begin{equation} \label{S} 
S=\{u\in H^1_0 (\Omega ) ,\int |u|^{\alpha +2}=1\}.
\end{equation} 
Using that $H^1 (\Omega  )$ is compactly injected in $L^{\alpha +2} (\Omega )$ one easily sees that \eqref{min} has a (unique) positive solution $\tilde U$. It is also clear that $U=k\tilde U$ solves \eqref{ellip} for a judicious choice of $k$.  It then follows from standard symmetrization arguments that $U$ is radial and radially decreasing. 

One may also obtain $U$ as a mountain pass solution. Consider 
\begin{equation} \label{E} 
E(u)=\frac {\rho } {2}\int  |u|^2+\frac {1} {2}\int  |\nabla u|^2-\frac {1} {\alpha +2}\int |u|^{\alpha +2}.
\end{equation} 
and 
\begin{equation} \label{Gamma} 
\Gamma =\{\gamma \in C([0,1]; H^1_0 (\Omega )), \gamma (0)=0, \gamma (1))=u_1\}, 
\end{equation}  
where $E(u_1)<0$. Then $E(u)$ is well-defined for $u\in H^1_0 (\Omega ) $ and $\Gamma $ is nonempty. In addition, 
\begin{equation} \label{c} 
c=\inf_{\gamma \in \Gamma } \sup_{t\in [0,1] } E(\gamma (t))
\end{equation} 
is a critical value of $E$ such that $c=E(U)>0$ and $E'(U)=0$. Moreover, it can be easily shown that $U$ is a ground state solution, i.e., $E(U)\le E(V)$ for all solution $V\ne 0$ of \eqref{ellip}.   

The general picture essentially remains unchanged for real solutions $u$ of \eqref{ellip} in the whole space $\Rn$, provided $\rho >0$ and $u\in H^1_{\rm rad} (\Rn )$, the space of radially symmetric functions of $H^1(\Rn)$. In fact, using that $H^1_{\rm rad} (\Rn )$ is compactly injected in $L^{\alpha +2} (\Rn )$ \cite{strauss}, the existence of a positive solution $U$ can be obtained either by solving \eqref{min} or \eqref{c}, when $S$ and $\Gamma $ are  redefined by replacing $H^1_0 (\Omega ) $ by $H^1_{\rm rad} (\Rn )$, see \eqref{S}, \eqref{Gamma}. Again, symmetrization arguments ensure that $U$ is radially decreasing. In fact, $U$ decays exponentially \cite{cazenave_s}. It is not difficult to see that both methods provide the same solution. However, in \cite{kwong} it is also shown the uniqueness of positive radially symmetric solutions in the whole space. (An alternative variational characterization of $U$ involving the so-called Gagliardo-Nirenberg quotient is presented in \cite[Proposition~2.6] {weinstein}. 

For $\Omega $ either the unitary ball or the whole space, we consider the linearized operator   
\begin{equation} \label{lin} 
Lv= \rho v-\Delta v  - U^{\alpha }v - \alpha U ^{\alpha} \Re v,
\end{equation}   
where $U$ is the positive solution of \eqref{ellip}.  More precisely, we set $D(L )=H^2(\Omega )\cap H^1_0(\Omega )$ and define $L :D(L )\subset L^2 (\Omega )  \to L^2 (\Omega )$ by \eqref{lin}.  Then $Lv=L_+\Re v+i L_-\Im v$ where  
\begin{equation} \label{L+}
L_+v=\rho v-\Delta v-(\alpha +1) |U  |^{\alpha} v,
\end{equation}  
\begin{equation}  \label{L-} 
L_-v=\rho v-\Delta v- |U  |^{\alpha} v.
\end{equation}
 We study below the operators $L_+$ and $L_-$. 
 \section{The linearized operator: the case of $\Omega =\Rn$}\label{Rn} 
In the case $\Omega=\Rn$, under a suitable rescaling we may assume that $\rho=1$ in \eqref{lin}.  
We want to show that $L_+$ is an injective operator when restricted to the space $V:=L_{\rm rad}^2(\Rn)$ of  radially symmetric and square integrable functions. We define $D(L_+)=   H^2 (\Rn )\cap V$ and consider $L_+ : D(L_+)\subset V\to V$ given by \eqref{L+}.  

Set $\sigma =U^\alpha $ and denote $V_\sigma $  the space $L^2 (\Rn )$ equipped with the scalar product
\begin{equation} 
\langle u,v\rangle_\sigma  =\int \sigma uv.
\end{equation} 
We also introduce $K: V_\sigma \to V_\sigma $ such that for $v\in V_\sigma$
\begin{equation} \label{K} 
K v=(\alpha +1)(I-\Delta )^{-1} U^\alpha v.
\end{equation} 
We have that $K$ is a positive, symmetric operator. Using that $U$ decays to zero at infinity, a standard argument shows that $K$ is compact. Denote $\{\varphi _j\}_{j\in \N}$ the orthonormal basis of eigenvectors of $K$ and $\{\lambda _j\}_{j\in \N}$ the corresponding set of eigenvalues.  Then $\lambda _j>0$ and $K\varphi _j=\lambda _j\varphi _j$ is equivalent to
\begin{equation} \label{equiv} 
\varphi _j-\Delta \varphi _j=\frac {\alpha +1} {\lambda _j}U^\alpha \varphi _j.
\end{equation} 
Therefore $\varphi _1=U$ and $\lambda _1=\alpha +1$. We will now prove  that $\lambda _2\le 1$. This is a consequence of the fact that $U$ is a mountain-pass solution. We present a simple proof below, which uses the specific form of the function $E(u)$. For the proof that general critical points of mountain-pass type have Morse index equal to one, see \cite{hofer}   
\begin{lem}\label{lambda_2} 
$\lambda _2\le 1$.
\end{lem}
\begin{proof} 
We first remark that for $k>1$ large enough $\gamma_0 (t)=ktU\in \Gamma $, see \eqref{Gamma}. In addition, 
$\max_{u\in \gamma_0} E(u)=E(U)$.   

We argue by contradiction and suppose that $\lambda _2> 1$. We get from \eqref{equiv} that
\begin{equation} 
\langle L_+\varphi _2,\varphi _2\rangle =(\alpha +1)(\lambda _2^{-1}-1)\int U^\alpha \varphi _2^2<0.
\end{equation}  
Consider now the plane $\pi$ containing $U$ and $\varphi _2$ and let 
\begin{equation} 
\gamma _1=\{U-\delta (\cos t U+\sin t \varphi _2), t\in [0,\pi]\}.
\end{equation} 
be an arc of circle in $\pi$ joining $(1-\delta )U$ and $(1+\delta )U$. We have that 
\begin{equation} \label{taylor} 
E(u)=E(U)+E'(U)(u-U)+\frac {1} {2}\langle L_+(u-U),u-U\rangle +o(\delta ^2).
\end{equation} 
Moreover, using that $\langle L_+U,\varphi _2\rangle=-\alpha \langle U,\varphi _2\rangle_\sigma =0$  we get
\begin{equation} 
\langle L_+(u-U),u-U\rangle =\delta ^2(\cos^2 t \langle L_+U,U\rangle+\sin^2 t \langle L_+\varphi _2,\varphi _2\rangle)<0.
\end{equation} 
Using this, \eqref{taylor} and that $E'(U)=0$, we see that we can choose $\delta $ small enough so that $E(u)<E(U)$ for all $u\in \gamma _1$.   

Let now $\gamma $ be the curve obtained by replacing the path of $\gamma _0$ going from $(1-\delta )U$ to $(1+\delta )U$ by $\gamma _1$. Then $\gamma \in \Gamma $ and $\max_{u\in \gamma } E(u)<E(U)$, leading to a contradiction. This shows that $\lambda _2\le 1$.
\end{proof} 
\begin{lem} \label{one_zero} 
Suppose $L_+\varphi =0$, $\varphi \ne 0$. Then there exists a unique $r^*>0$ such that $\varphi (r^*)=0$. 
\end{lem} 
\begin{proof} 
Let $B_R$ be the ball of radius $R$ of $\Rn$. For $v\in L^2_{\rm rad}(B_R)$ let $u\in H^2_{\rm rad}(B_R)\cap H^1_0 (B_R )$ satisfy $(I-\Delta )u=(\alpha +1)U^\alpha v$. We define $K_R: L^2_{\rm rad}(B_R)\to L^2_{\rm rad}(B_R)$ such that $K^R v=u$. Then $K^R$ is a compact operator, which is symmetric and positive for the scalar product
\begin{equation} 
\langle u,v\rangle _{\sigma ,R}=\int_{B_R} U^\alpha uv.
\end{equation} 
Denote  $\{\varphi _j^R\}_{j\in \N}$ an orthonormal basis of eigenvectors of $K^R$, associate to the set $\{\lambda _j^R\}_{j\in \N}$ of eigenvalues, so that  
\begin{equation} \label{util} 
(I-\Delta )\varphi _j^R=\frac {\alpha +1} {\lambda _j^R} U^\alpha  \varphi _j^R
\end{equation} 
and 
\begin{equation} \label{normal} 
\int_{B_R} U^\alpha \varphi _i^R\varphi _j^R=\delta _{ij}.
\end{equation}  
Moreover, it is easy to see that 
\begin{equation} \label{lambda_R} 
\lambda _j^R<\lambda _1^R<\lambda _1,
\end{equation} 
where $\lambda _1$ is the first eigenvalues of $K$ given by \eqref{K}, and that for all $j$ 
\begin{equation} \label{conv} 
\lambda _j^R\nearrow \lambda _j^\infty\le \lambda _1
\end{equation}     
as $R\to \infty$. We extend $\varphi _j^R(r)=0$ for $r>R$. Using \eqref{util} and \eqref{normal} we see that 
\begin{equation} 
\int | \varphi _j^R|^2+|\nabla  \varphi _j^R|^2=\frac {\alpha +1} {\lambda _j^R}\int_{B_R}U^\alpha |\varphi _j^R|^2=\frac {\alpha +1} {\lambda _j^R}.
\end{equation} 
It follows then from \eqref{conv} that $\{\varphi _j^R\}_{R\ge  \underline{R}}$ is uniformly bounded in $H^1 (\Rn )$ for all $\underline{R}>0$. Upon considering a subsequence, we may write that there exists $\varphi^\infty $ in  $H^1 (\Rn )$ such that $\varphi _j^R\rightharpoonup \varphi_j^\infty $ weakly in  $H^1 (\Rn )$ as $R\to \infty$. Using that $U(r)\to 0$ as $r\to \infty$, we readily obtain that 
\begin{equation} \label{} 
\int_{B_R} U^\alpha \varphi _i^\infty\varphi _j^\infty=\delta _{ij},
\end{equation}  
with
\begin{equation} \label{} 
(I-\Delta )\varphi _j^\infty=\frac {\alpha +1} {\lambda _j^\infty} U^\alpha  \varphi _j^\infty.
\end{equation} 
Thus, $\lambda _j^\infty$ is an eigenvalue of $K$, associated to $\varphi  _j^\infty$. 
 
Suppose now that $L_+\varphi =0$ so that $K\varphi =\varphi $. Assume that $\varphi (\rho )=\varphi (R)=0$ for some $0<\rho <R$.  Then $K^R\varphi =\varphi $. Thus $1=\lambda^R _j$ for some $j$. However, since $\varphi $ changes sign once in $(0,R)$,   $1=\lambda^R _2$. Thus $\lambda_2 ^\infty>1$. Since $\lambda_2 ^\infty\ne \lambda _1$, we see that 
$1<\lambda_2 ^\infty\le \lambda _2$. But this contradicts Lemma~\ref{lambda_2}. Thus, $\varphi (r)$ has a single zero $r^*>0$.
\end{proof} 

We next present the ingenious argument of \cite{cgnt} to show that $L_+$ is injective.  
\begin{lem}\label{ker_L+} 
$L_+$ is injective.
\end{lem}
\begin{proof} 
We argue by contradiction and assume that there exists $\varphi\ne 0$ such that $L_+\varphi=0$. Using Lemma~\ref{one_zero}, we may assume that there exists $r^*>0$ such that $\varphi (r)>0$ for $r<r^*$ and  $\varphi (r)<0$ for $r>r^*$. Set now
\begin{equation} \label{eta} 
\eta (x)=U(x)+\frac {\alpha } {2} x. \nabla U(x)
\end{equation} 
Since $U$ decays exponentially, $\eta \in H^1_{\rm rad} (\Rn )$. Moreover, a straightforward calculation yields 
\begin{equation} \label{L_eta} 
L_+\eta   = - \alpha U.
\end{equation} 
Define $w=U^{\alpha }(r^*)\eta-U$ and $z=L_+ w$. Then $z(r)= \alpha U(r)(U^{\alpha }(r)-U^{\alpha }(r^*)) $ and so $z (r)>0$ for $r<r^*$ and  $z (r)<0$ for $r>r^*$. Hence, $z(r)\varphi (r)>0$ for $r\ne r^*$. However,  this is in contradiction with the fact that 
\begin{equation*} 
\langle \varphi_2, z \rangle= \langle \varphi_2,L_+ w\rangle=\langle L_+\varphi_2, w\rangle=0.
\end{equation*} 
This shows that $L_+$ is injective.
\end{proof} 
Using decomposition in spherical harmonics, in \cite{weinstein} and in \cite{cgnt} it is proved that the complete kernel of  $L_+$ in $L^2(\Rn)$ is $\ker L_+=[\partial _1 U,\partial _2U,\cdot,\partial _NU]$. Note that $\partial _jU$ is not a radial function.   

We may now characterize the kernel of $L$ given by \eqref{lin}.  
\begin{prop}\label{ker_L_Rn} 
We have $\ker L=[iU]$.
\end{prop}
\begin{proof} 
If $v\in \ker L$ then $\Re v\in \ker L_+$ and $\Im v\in \ker L_-$. It follows from Lemma~\ref{ker_L+} that $\Re v=0$. Moreover, if $\varphi \in \ker L_-$ then $\varphi $ is an eigenvalue of $K$ given by \eqref{K}, associated to $\lambda =\alpha +1$. But  $\alpha +1$ is the first eigenvalue of $K$ and $KU=(\alpha +1)U$. Hence $\ker L_-=[U]$ so that  $\ker L=[iU]$.
\end{proof} 

\section{The linearized operator: the case of a ball}\label{ball} 
Let $\Omega\subset \Rn$ be the unitary ball and suppose \eqref{rho2} holds.  Let $U$ be the unique positive solution of \eqref{ellip} and let $L$ be given by \eqref{lin}. Then $v\in \ker L$ if and only if $\Re v\in \ker L_+$ and $\Im v\in \ker L_-$. Since $L_-U=0$ and $U>0$, it follows that $\ker L_-=[U]$ is a one-dimensional subspace. We will now show that $L_+$ is injective. This is proved in \cite{dgp} for $\rho=0$. For the reader's convenience, we reproduce the arguments here. The two preliminary results, Lemma~\ref{eLem0} and Lemma~\ref{eLem0} hold in fact for $\rho > -\lambda _1$ and will be useful in the proof of the general case.

For  $x\in \R^N$ write $x=(t,y)$, where $t\in \R$, $y\in \R^{N-1}$, if $N>1$ or $x=t$ if $N=1$. Set $\Omega ^*=\{x\in \Omega ,t<0\}$, $D(L_+^*)=H^2(\Omega ^*)\cap H^1_0(\Omega ^*)$ and $L_+^*:D(L_+^*)\subset L^2 (\Omega ^*)  \to L^2 (\Omega ^*)$ be given by \eqref{L+}.   
\begin{lem} \label{eLem0} 
We have $\lambda _1^*=\lambda _1(L_+^*)>0$.
\end{lem}   
\begin{proof} 
 Let $v=\partial _t U$. It is well known that $v>0$ over $\Omega ^*$ with $v>0$ over $\Gamma ^*=\{x\in  \overline{\Omega ^*} , |x|=1\}$. Moreover, taking the derivative with respect to $t$ in \eqref{ellip}  we see that $L_+^* v=0$. Consider $u_1$ a positive eigenvector of $L_+^*$, so that $L_+^*u_1=\lambda _1^* u_1$. Then
 \begin{equation} 
\lambda _1^* \int _{ \Omega ^* } u_1v=  \int _{ \Omega ^* } v L_+^*u_1=- \int _{\partial  \Omega ^* } v \partial _\eta u_1>0.
 \end{equation} 
 This shows that $\lambda _1^*>0$. 
\end{proof} 
As a consequence, we have the following.
\begin{lem} \label{eLem00} 
Let $v$ satisfy $L_+v=0$. Then $v$ is radially symmetric.
\end{lem}
\begin{proof} 
If $v\in \ker L_+$ then $v\circ R\in \ker L_+$ for all unitary transformation $R$. It thus suffices to show that $v(t,y)$ is symmetric with respect to $t$. 
Define $\psi (x)=v(t,y)-v(-t,y)$. Then $L_+^*\psi=0$, with $\psi=0$ over $\partial \Omega ^*$. It follows from Lemma~\ref{eLem0} that $\psi=0$. This ends the proof.
\end{proof} 
\begin{lem} \label{eLem1} 
Suppose $\rho=0$. Then the operator $L_+$ given by \eqref{L-} is injective.  
\end{lem}
\begin{proof} 
Let $v\in H^2(\Omega )\cap H^1_0(\Omega )$ satisfy $L_+v=0$. Then 
\begin{equation} \label{fLem1.1} 
0=\int_\Omega L_+v\, U=\int_\Omega L_+U\, v=-\alpha \int_\Omega  U^{\alpha +1}v.
\end{equation} 
Consider now the Pohozaev function $\psi=x.\nabla U$. We have that 
$\partial _j\psi =\partial _j U+x.\nabla \partial _j U$, so that 
$\partial^2 _{jj}\psi=2\partial^2 _{jj}U+x.\nabla \partial^2 _{jj}U$. Thus, 
$\Delta \psi =2\Delta U+x.\Delta \nabla U$. In addition, we get from \eqref{ellip} that
$\Delta \nabla U+ (\alpha +1)U^\alpha\nabla U=0$, and so 
\begin{equation} \label{fLem1.2} 
\Delta \psi =2\Delta U-(\alpha +1)U^\alpha x.\nabla U=-2U^{\alpha +1}-(\alpha +1)U^\alpha \psi .
\end{equation}  
It follows from \eqref{fLem1.2} and \eqref{fLem1.1} that
\begin{equation} 
0=2\int_\Omega  U^{\alpha +1}v=\int_\Omega L_+\psi \, v=
\int_{\partial \Omega }\psi \,\partial _\eta v=
 \int_{\partial \Omega }\partial _\eta U \,\partial _\eta v. 
\end{equation}
We conclude from Lemma~\ref{eLem00} and Hopf's strong maximum principle that $v'(1)=0$. Therefore, $v=0$. 
\end{proof} 
\begin{rem} \label{rem-starsh}
Following \cite{dgp}, Lemma~\ref{eLem1} holds in the case $N=2$ and $\Omega$ any regular bounded star shaped domain.
\end{rem}
The proof that $L_+$ injective in the case $\rho\ne 0$ follows the arguments of  \cite{cgnt}, see Lemma~\ref{ker_L+}. 
\begin{lem}\label{eLem2} 
Assume $\rho>-\lambda _1$, $\rho \ne 0$. Then $L_+$ is injective.
\end{lem}
\begin{proof}
Let   
\begin{equation} \label{eta} 
\eta (r)=2U(r)+\alpha rU'(r). 
\end{equation} 
A straightforward calculation gives that
\begin{equation} \label{L_eta} 
L_+\eta   = -2\rho \alpha U.
\end{equation} 
We want to show that $\ker{L_+}=\{0\}$. We argue by contradiction and assume that there exists $\varphi_2\ne 0$ such that $L_+\varphi_2=0$. Since $U$ is a mountain-pass solution of~\eqref{ellip}, we know that $\lambda_2(L_+)\ge 0$, see  \cite{hofer}. Thus, $\varphi _2$ is an eigenvector associated to the second eigenvalue $\lambda _2=0$. By Lemma~\ref{eLem00}, $\varphi_2$ is radial. Using standard comparison arguments, it is easy to see that $\varphi_2$ has a single zero $r_0$ in $(0,1)$.  For $b\in (0,1)$, define 
\begin{equation} 
g(r)=
\begin{cases}
1 &  0\le r<b,\\
2\frac {1-r} {1-b}-\frac {(1-r)^2} {(1-b)^2} & b\le r\le 1.
\end{cases} 
\end{equation} 
Then $g''(r)=g'(r)=0$ if $r<b$. For $r>b$,
\begin{equation} \label{g'} 
g'(r)= \frac {2(b-r)} {(1-b)^2},
\end{equation} 
\begin{equation} \label{g''} 
g''(r)=-\frac {2} {(1-b)^2}.
\end{equation} 
We remark that we can choose $b$ close enough to $1$ so that
\begin{equation} 
b>r_0
\end{equation} 
and 
\begin{equation} \label{b} 
U^\alpha (r)<U^\alpha (r_0)g(r) 
\end{equation} 
for $r>b$.
Set now $w(r)=g(r)\eta (r)$, see \eqref{eta}. It follows from \eqref{L_eta} that 
\begin{equation} \label{L_w} 
L_+w=gL_+\eta -2\nabla g.\nabla \eta -\eta\Delta g= -2g\rho \alpha U-2\nabla g.\nabla \eta -\eta\Delta g.
\end{equation}  
Thus $L_+w=-2\rho \alpha U$ if $|x|<b$. 

For $|x|>b$, we get from \eqref{eta}, \eqref{g'} and \eqref{g'} that  
\begin{equation} \label{first} 
\nabla g(r).\nabla \eta(r)=g'(r)\eta'(r)= \frac {2(b-r)} {(1-b)^2}((2+\alpha )U'(r)+\alpha r U''(r) )
\end{equation} 
and that
\begin{equation} \label{second} 
\Delta g=g''+\frac {N-1} {r}g'= -\frac {2} {(1-b)^2}-\frac {2(N-1)(r-b)} {r(1-b)^2}.
\end{equation} 
Defining $h=-2\nabla g.\nabla \eta -\eta\Delta g$ we get from \eqref{eta}, \eqref{first} and \eqref{second} that there exists $K>0$ such that 
\[
(1-b)^2h(r)\le 4U(r)+2\alpha rU'(r)+K(r-b).
\]     
Since $U(1)=0$ and $U'(1)<0$, $1-b$ can be taken eventually smaller so that  
\begin{equation} \label{hnegative} 
h(r)<0 \, \, \hbox{for} \, \, r>b.
\end{equation}     
Set now $t  =U^\alpha (r_0)/(2\rho )$ and $z=L_+(-U+t w)$. Hence,  by  \eqref{L_w} we get  
\begin{equation} 
z=\alpha U^{\alpha +1}+t(-2g\rho \alpha U +h).
\end{equation} 
Let us show that $z$ and $\varphi _2$ have the same sign. 
For $r<b$ we use that $g=1$ and $h=0$  to get    
\begin{equation*} 
z(r)=\alpha U(r)(U^{\alpha }(r)-U^\alpha (r_0)). 
\end{equation*} 
It follows that $z(r)>0$ if $r<r_0$ and $z(r)<0$ for $r\in (r_0,b)$. In addition, using \eqref{hnegative} and \eqref{b}, we get for $r>b$ that   
\begin{equation*} 
z(r)<\alpha U(r)(U^{\alpha }(r)-gU^\alpha (r_0))<0. 
\end{equation*} 
We see then that $z(r)\varphi_2(r)>0$ for $r\ne r_0$. But   
\begin{equation*} 
\langle \varphi_2,z\rangle =\langle \varphi_2,L_+(-U+\beta w)\rangle=\langle L_+\varphi_2,-U+\beta w\rangle=0, 
\end{equation*} 
giving a contradiction. This shows that $L_+$ is injective.
\end{proof}
We present now the main result of this section.
\begin{prop}\label{ker_L_ball} 
We have $\ker L=[iU]$.
\end{prop}
\begin{proof} 
Let $v\in \ker L$. Then $\Re v\in \ker L_+$ and $\Im v\in \ker L_-$. It follows from Lemma~\ref{eLem1} and Lemma~\ref{eLem2} that $\Re v=0$. Moreover, as discussed in the beginning of this section $\ker L_-=[U]$. This closes the proof.
\end{proof} 
\section{Proofs of Theorem~\ref{thm1} and Theorem~\ref{thm2} }\label{proofs} 
In this section we denote $L^p (\Omega )$ the real Banach space whose elements are complex-valued  functions. In particular, $L^2 (\Omega ) $ is a Hilbert space for the scalar product 
\begin{equation} \label{sp} 
(u,v)=\Re \int_\Omega u \overline{v}.
\end{equation} 
Accordingly, $H^m(\Omega )$ denotes a real Hilbert space having complex elements. 
\begin{proof} [Proof of  Theorem~\ref{thm1}]
For a fixed $\theta \in (-\pi /2,\pi /2)$ set $X=\R\times (H_{\rm rad}^2(\Rn ))$ and $F:(-\pi /2, \pi/2) \times X\to L_{\rm rad}^2(\Rn )$ such that 
\begin{equation} \label{fF2} 
F(\gamma , \omega , u) =\rho u- \Delta u  - e^{i(\gamma -\theta )}  |u |^{\alpha }u  -i \omega e^{ -i \theta } u.
\end{equation} 
Note that $F$ is well defined due to Sobolev embedding $H^2(\Rn )\hookrightarrow L^{ 2\alpha +2}(\Rn )$.  

Then $\varphi _\gamma =e^{i\omega _\gamma t}u_\gamma $ is a solution of \eqref{CGL} if and only if $F(\gamma , \omega _\gamma , u_\gamma )=0$.  Note that $F(\theta ,0,Ue^{i\beta })=F(\theta ,g(\theta ))=0$. In addition, it is immediate to see that $F$ is a $C^1$ function such that
\begin{align*}
\frac {\partial F} {\partial \omega  } (\gamma , \omega , u)\mu  &=  
-i e^{ -i \theta } u \mu , \\
\frac {\partial F} {\partial u }(\gamma  , \omega , u) v&= 
\rho v-\Delta v  -e^{i(\gamma -\theta )} \bigl[   |u |^{\alpha }v 
+   \alpha    |u |^{\alpha -2} u\Re ( \overline{u}  v)\bigr] -i \omega e^{-i\theta }v . 
\end{align*}  
By the surjective form of the Implicit Function Theorem~\cite[Theorem~4.H, p.177] {Zeidler}, the proof will be completed once we show that 
$\partial _{\omega ,u} F(\theta ,0,Ue^{i\beta}): X\to L^2 (\Omega ) $ is surjective. Note that  
\begin{align}
\frac {\partial F} {\partial \omega  } (\theta  , 0 , Ue^{i\beta }) &=  
-i e^{ -i \theta } Ue^{i\beta },\label{fthm2} \\
\frac {\partial F} {\partial u }(\theta   , 0 , Ue^{i\beta }) v&= 
\rho v-\Delta v  - U ^{\alpha }v 
-   \alpha    U ^{\alpha}e^{i\beta } \Re (e^{-i\beta }v) ,\label{fthm3} 
\end{align}  
so that 
\[\partial_{ \omega,u} F (\theta  , 0,Ue^{i\beta })(\mu,v)=e^{i\beta }\partial_{ \omega,u} F (\theta  , 0,U)(\mu,e^{-i\beta }v).
\]
It thus suffices to consider the case $\beta =0$. 

Given $f\in L_{\rm rad}^2 (\Omega ) $, $\partial_{ \omega,u} F (\theta  , 0,U)(\mu,v)=f$ is equivalent to   
 \begin{equation} \label{bijective} 
 -i e^{ -i \theta } U\,\mu+L v =f,
 \end{equation}  
where $L$ is given by \eqref{lin}.  Note that $L$ is a self-adjoint operator in $L_{\rm rad}^2 (\Rn ) $ for the scalar product \eqref{sp}. Using that $\ker L=[iU]$, see Proposition~\ref{ker_L_Rn}, we choose $\mu$ such that  
\[
\tilde f=f+i e^{ -i \theta } U\,\mu\in (iU)^\perp,
\] 
i.e., 
\begin{equation} \label{mu} 
\mu=-\frac {1} {\cos\theta \|U\|^2_{ L^2 (\Rn )  }}\int _{ \Rn  } \Im fU.
\end{equation} 
The fact that $Lv=\tilde f$ has a solution for $\tilde f\in (iU)^\perp$ follows from the Fredholm Alternative applied to the compact operator $K=(\rho -\Delta )^{-1}U^\alpha $, see Section~\ref{Rn}. This shows that $L$ is surjective and closes the proof.
\end{proof} 
\begin{proof} [Proof of  Theorem~\ref{thm2}]
Set $X=\R\times (H^2(\Omega )\cap H^1_0(\Omega ))$ and 
define $F:(-\pi /2, \pi/2) \times X\to L^2(\Omega )$ by \eqref{fF2}. 
The arguments of the proof of Theorem~\ref{thm1} are still valid in this case where $\Omega $ is a ball. The fact that  $\ker L=[iU]$ was stablished in Proposition~\ref{ker_L_ball}. 
\end{proof} 
\begin{rem} \label{rem}
%The curve $g_\gamma =(\omega _\gamma ,u _\gamma)$ in Theorem~\ref{thm2} is not unique. Indeed, we may take $w\in L^2 (\Omega ) $ such that $\langle iU,w\rangle \ne 0$ and consider $Y_w=H^2(\Omega )\cap H^1_0(\Omega )\cap  w^\perp$. Given $\tilde f\in (iU)^\perp$, let $\tilde z\in (iU)^\perp$ be the unique solution of $L\tilde z=\tilde f$. Setting $z=\tilde z-tiU$, where 
%\[
%t=\frac {\langle \tilde z,w\rangle} {\langle iU,w\rangle},
%\]
%we see that $z\in w^\perp$. We may then modify the proof of Theorem~\ref{thm2} by applying the standard Implicit Function Theorem to find a curve of standard 

\begin{enumerate}

\item{} Let $\Omega $ be the unitary ball of $\Rn$ and let $\tilde X=\R\times H^2(\Omega )\cap H^1_0(\Omega )\cap  (iU)^\perp$. Given $\tilde f\in (iU)^\perp$, there exists a unique $\tilde z\in (iU)^\perp$ such that $L\tilde z=\tilde f$. We may thus    modify the proof of Theorem~\ref{thm2} and apply the standard Implicit Function Theorem to find a unique curve of 
$g_\gamma =(\omega _\gamma ,u _\gamma)$ in $\tilde X$ such that $\varphi _\gamma =e^{i\omega _\gamma} u _\gamma$ is a standing wave solution of \eqref{CGL}.  Since the equation is invariant under unitary transformations and $U$ is radially symmetric, it follows by the uniqueness of $g_\gamma $ that 
$\varphi _\gamma$ is radially symmetric. 
\item{} Theorem~\ref{thm2} is still valid in the case $N=2$, $\Omega$ any bounded regular star shaped domain and $\rho =0$. From Remark \ref{rem-starsh} the argument of the proof applies without any change.
\end{enumerate}
\end{rem}

\end{document}